\documentclass[12pt]{amsart}
\usepackage[all]{xy}
\usepackage{amscd,amssymb,amsmath}
\newtheorem{theorem}{Theorem}

\newtheorem{que}[theorem]{Question}

\theoremstyle{definition}

\newtheorem{rk}[theorem]{Remark}
\newtheorem{ex}[theorem]{Example}

\theoremstyle{definition}

\def\<{\langle}
\def\>{\rangle}

\def\b{\beta}
\def\g{\gamma}

\def\r{\rho}

\def\s{\sigma}

\def\G{{\Gamma}}

\renewcommand{\to}{\rightarrow}

\DeclareMathOperator{\re}{re}

\def \bd {\partial}
\def \hd {Heegaard decomposition}
\def \re {representation}
\def \CR {\mathcal R}

\def\Fix{\operatorname{Fix}}

\newcommand{\bea} {\begin{eqnarray*}}
\newcommand{\beq} {\begin{equation}}
\newcommand{\bey} {\begin{eqnarray}}
\newcommand{\eea} {\end{eqnarray*}}
\newcommand{\eeq} {\end{equation}}
\newcommand{\eey} {\end{eqnarray}}

\begin{document}

\title[Nielsen number is a  knot  invariant]
{Nielsen number is a knot  invariant  }

\author{Alexander Fel'shtyn}

\address{Instytut Matematyki, Uniwersytet Szczeci\'nski ul. Wielkopolska 15,
70-451 Szczecin, Poland and Boise State University, 1910
University Drive, Boise, Idaho, 83725-155, USA  }
\email{felshtyn@diamond.boisestate.edu; felshtyn@mpim-bonn.mpg.de}

\begin{abstract}
  We show that the  Nielsen number is a  knot  invariant  via representation variety
\end{abstract}
\maketitle

\section{Introduction}
We briefly describe the few
basic notions of Nielsen fixed point theory(see \cite{f1}).
We assume  $X$ to be a connected, compact
polyhedron and $f:X\rightarrow X$ to be a continuous map.
Let $p:\tilde{X}\rightarrow X$ be the universal cover of $X$
and $\tilde{f}:\tilde{X}\rightarrow \tilde{X}$ a lifting
of $f$, i.e. $p\circ\tilde{f}=f\circ p$.
Two liftings $\tilde{f}$ and $\tilde{f}^\prime$ are called
{\sl conjugate} if there is a $\gamma\in\Gamma\cong\pi_1(X)$
such that $\tilde{f}^\prime = \gamma\circ\tilde{f}\circ\gamma^{-1}$.
The subset $p(Fix(\tilde{f}))\subset Fix(f)$ is called
{\sl the fixed point class of $f$ determined by the lifting class $[\tilde{f}]$}.Two fixed points $x_0$ and $x_1$ of $f$ belong to the same fixed point class iff  there is a path $c$ from $x_0$ to $x_1$ such that $c \cong f\circ c $ (homotopy relative endpoints). This fact can be considered as an equivalent definition of a non-empty fixed point class.
 Every map $f$  has only finitely many non-empty fixed point classes, each a compact  subset of $X$.
A fixed point class is called {\sl essential} if its index is nonzero.
The number of essential fixed point classes is called the {\sl Nielsen number}
of $f$, denoted by $N(f)$.The Nielsen number is always finite.  $N(f)$ is a  homotopy invariant.
In the category of compact, connected polyhedra, the Nielsen number
of a map is, apart from certain exceptional cases,
 equal to the least number of fixed points
 of maps with the same homotopy type as $f$.\\
Let us consider a braid representative of a knot and induced map of
corresponding representation variety(see section 2).We prove in section 3
that  the Nielsen number of induced map is a invariant under Markov moves and
so is a knot invariant.\\
The author came to the  idea that Nielsen number is a  knot  invariant  at the summer  2003, after conversations with Jochen Kroll and Uwe Kaiser.
The author would like to thank the Max-Planck-Institute f\"ur Mathematik, Bonn
for kind hospitality and support.

\section{ Casson-Lin invariant}

We recall firstly the  Lin's construction in \cite{lin}
for the
intersection number of the representation spaces corresponding to a braid
representative of a knot $K$ in $S^3$. 
Let $(S^3, D^3_+, D^3_-, S^2)$ be a Heegaard decomposition of $S^3$ with genus
$0$, where
\[ S^3 = D^3_+ \cup_{S^2} D^3_- , \ \ \ \bd D^3_+ = \bd D^3_- =
D^3_+ \cap D^3_- = S^2 .\]
Suppose that a knot  $K \subset S^3$ is
in general position with respect to this Heegaard decomposition. So
$K \cap S^2 = \{x_1, \cdots, x_n, y_1, \cdots, y_n\}$, $K \cap D^3_{\pm }$
is a collection of unknotted, unlinked arcs
$\{ \g_1^{\pm }, \cdots, \g_n^{\pm } \} \subset D^3_{\pm }$,
where $\bd \g_i^- = \{x_i, y_i\}$ and $\{\g_1^+, \cdots, \g_n^+\} = K \cap
D^3_+$ becomes a braid of $n$ strands inside $D^3_+$. Denote by  $\b $ a corresponding word  in the
braid group $B_n$. For the top end points  $x_i$ of $\g_i^+$, the bottom
end points of $\{\g_1^+, \cdots, \g_n^+\}$ give  a permutation of
$\{y_1, \cdots , y_n\}$ which generates a map
\[ \pi : B_n \to S_n ,\]
where $\pi (\b )$ is the permutation of $\{y_1, \cdots , y_n\}$ in the symmetric
group of $n$ letters. Let $K = \overline{\b }$ be the closure of $\b$.
It is well-known that there is a correspondence between a knot and a braid $\b$
with $\pi (\b )$ is a complete cycle of the $n$ letters (see \cite{bi}).

There is a corresponding {\hd} for the complement of a  $K$,
\[ S^3 \setminus K = (D^3_+ \setminus K) \cup_{(S^2 \setminus K)} (D^3_-
\setminus K ) , \]
\[  D^3_{\pm } \setminus K = D^3_{\pm } \setminus (D^3_{\pm }
\cap K), \ \ \ S^2 \setminus K = S^2 \setminus (S^2 \cap K) . \]
Thus by  Seifert-van Kampen theorem we have following diagramm
\[ \begin{array} {ccc}
\pi_1 (S^2 \setminus K) & \to & \pi_1 (D^3_+ \setminus K)\\
\downarrow & & \downarrow \\
\pi_1 (D^3_- \setminus K ) & \to & \pi_1 (S^3 \setminus K) ,
\end{array} \]
and a corresponding diagramm  of {\re} spaces
\begin{equation} \label{pullback}
 \begin{array} {ccc}
{\CR}(S^2 \setminus K)& \leftarrow & {\CR}(D^3_+ \setminus K)\\
\uparrow & & \uparrow \\
{\CR}(D^3_- \setminus K ) & \leftarrow & {\CR} (S^3 \setminus K) ,
\end{array} \end{equation}
where ${\CR }(X) = Hom (\pi_1 (X), SU(2))/SU(2)$ for $X =
S^2 \setminus K, D^3_{\pm } \setminus K, S^3 \setminus K$.

In \cite{ma}, Magnus used the trace free matrices to represent the
generators of a free group to show that the faithfulness of a
{\re} of braid groups in the automorphism groups of the rings generated
by the character functions on free groups. This is original idea
to have {\re}s with trace free along all meridians which Lin worked in
\cite{lin} to define the knot invariant. It has been carried out by
M. Heusener and J. Kroll  in \cite{hk} for the {\re} of knot groups with
the trace of the meridian fixed (not necessary zero).
Let ${\CR}(S^2 \setminus K)^{[i]}$ be the
space of $SU(2)$ {\re}s $\r : \pi_1 (S^2 \setminus K) \to SU(2)$
such that
\begin{equation} \label{trace}
 \r ([m_{x_i}]) \sim \left( \begin{array}{cc}
i & 0 \\ 0 & -i \end{array} \right) , \ \ \
\r ([m_{y_i}]) \sim \left( \begin{array}{cc}
i & 0 \\ 0 & -i \end{array} \right) , \end{equation}
where $m_{x_i}, m_{y_i}, i =1, 2, \cdots, n$ are the meridian circles around
$x_i, y_i$ respectively. Note that $\pi_1(S^2 \setminus K)$ is generated by
$m_{x_i}, m_{y_i}, i =1, 2, \cdots, n$
and one relation $\prod^n_{i=1} m_{x_i} = \prod^n_{i=1} m_{y_i}$.
Corresponding to (\ref{pullback}), we have
\begin{equation} \label{irre}
 \begin{array} {ccc}
{\CR}(S^2 \setminus K)^{[i]} & \leftarrow & {\CR}(D^3_+ \setminus K)^{[i]}\\
\uparrow & & \uparrow \\
{\CR}(D^3_- \setminus K )^{[i]} & \leftarrow & {\CR} (S^3 \setminus K)^{[i]} .
\end{array}
\end{equation}
The conjugacy class in $SU(2)$ is completely determined by its trace.
So the condition (\ref{trace}) can be reformulated for $\r \in {\CR}(X)^{[i]}$,
\begin{equation} \label{22}
trace \r ([m_{x_i}]) = trace \r ([m_{y_i}]) = 0 . \end{equation}
The space ${\CR}(S^2 \setminus K)^{[i]}$ can be identified with the space
of $2n$ matrices $X_1 \cdots, X_n, Y_1, \cdots, Y_n$ in $SU(2)$ satisfying
\begin{equation} \label{zero}
trace (X_i) = trace (Y_i) = 0, \ \ \ \ \mbox{for $i =1, \cdots, n$},
\end{equation}
\begin{equation} \label{product}
X_1 \cdot X_2 \cdots X_n = Y_1 \cdot Y_2 \cdots Y_n .
\end{equation}
Let $Q_n$ be the space $\{(X_1, \cdots , X_n) \in SU(2)^n | \ \
trace (X_i) = 0, i =1, \cdots, n \}$.
Let  ${\CR}^*(S^2 \setminus K)^{[i]}$ be the subset of  ${\CR}(S^2 \setminus K)^{[i]}$ consisting of irreducible representations.
Note that ${\CR}^*(S^2 \setminus K)^{[i]} = (H_n \setminus S_n)/SU(2)$
in Lin's notation \cite{lin}, where
\[H_n = \{(X_1, \cdots, X_n, Y_1, \cdots, Y_n) \in Q_n \times Q_n | \ \
X_1 \cdots X_n = Y_1 \cdots Y_n \},\]
$S_n$ is the subspace of $H_n$ consisting of all the reducible points.
Here $H_n \setminus S_n$ is the total space of a $SU(2)$-fiber bundle over
${\CR}^*(S^2 \setminus K)^{[i]}$.

Given $\b \in B_n$, we denote by $\G_{\b}$ the graph of $\b $ in $Q_n \times Q_n$,
i.e.
\[\G_{\b} = \{(X_1, \cdots, X_n, \b (X_1) , \cdots, \b (X_n)) \in Q_n \times Q_n\}.
\]
As an automorphism of the free group $Z[m_{x_1}]*Z[m_{x_2}]* \cdots *Z[m_{x_n}]$, this
element $\b \in B_n$ preserves the word $[m_{x_1}] \cdots [m_{x_n}]$. Thus we have
$ X_1 \cdots X_n = \b (X_1) \cdots \b (X_n), $
or in other words $\G_{\b}$ is a subspace of $H_n$. In fact, for $\overline{\b}
= K$, this subspace $\G_{\b}$ coincides with  the subspace of {\re}s
$\r : \pi_1(S^2 \setminus K) \to SU(2)$ in $H_n$ which can be extended
to $\pi_1(D^3_+ \setminus K)$,
$\G_{\b} = Hom (\pi_1(D^3_+ \setminus K), SU(2))^{[i]}.$
Hence the space ${\CR}^*(D^3_+ \setminus K)^{[i]} = \G_{\b, irre}/SU(2)$
is the irreducible $SU(2)$ {\re}s with traceless condition over
$D^3_+ \setminus K$.

In the special case $\b = id$, then $\G_{id}$ represents the diagonal in
$Q_n \times Q_n$,
\[\G_{id} = \{ (X_1, \cdots, X_n, X_1, \cdots, X_n) \in Q_n \times Q_n\} . \]
Since $K \cap D^3_-$ represents the trivial braid, this space $\G_{id}
\subset H_n$ can be identified with the subspace of {\re}s
in $Hom(\pi_1(S^2 \setminus K), SU(2))^{[i]}$ which  can be extended to
$\pi_1(D^3_- \setminus K)$, i.e.
$ \G_{id} = Hom (\pi_1(D^3_- \setminus K), SU(2))^{[i]} . $
By Seifert, Van-Kampen Theorem, the intersection
$\G_{\b} \cap \G_{id}$ is the same as the space of {\re}s of
$\pi_1(S^3 \setminus K)$ satisfying the monodromy condition $[i]$
(see (\ref{pullback})),
\[\G_{\b} \cap \G_{id} = Hom (\pi_1(S^3 \setminus K), SU(2))^{[i]} . \]

Given $\b \in B_n$ with $\overline{\b } = K$, there is an induced diffeomorphism
(still denoted by $\b $) from $Q_n$ to itself. Such a diffeomorphism also
induces a diffeomorphism $f_{\b }:  {\CR}^*(S^2 \setminus K)^{[i]} \to
{\CR}^*(S^2 \setminus K)^{[i]}$of the representation variety.

Note that
$\overline{\G}_{\b} = ({\G}_{\b} \setminus ({\G}_{\b} \cap S_n))/SU(2)$
is the image  of the  ``diagonal'' $\overline{\G}_{id}$ under diffeomorphism  $f_{\b }$.
By Seifert- van Kampen theorem (\ref{irre}), it is clear that the fixed point
set of $f_{\b }$ is
\[\ Fix (f_{\b }|_{{\CR}^*(S^2 \setminus K)^{[i]}}) =
\overline{\G}_{\b} \cap \overline{\G}_{id} =
{\CR}^*(S^3 \setminus K)^{[i]} .\]
 
The oriented submanifolds
$\overline{\G}_{\b}={\CR}^*(D^3_+ \setminus K)^{[i]},  \overline{\G}_{id}={\CR}^*(D^3_-\setminus K)^{[i]}$
intersects each other in a compact subspace of
${\CR}^*(S^2 \setminus K)^{[i]}$ from Lemma 1.6 in \cite{lin}.
Hence we can perturb $f_{\b }$, or
in another words perturb ${\CR}^*(D^3_+ \setminus K)^{[i]}$ to
$\hat{{\CR}}^*(D^3_+ \setminus K)^{[i]}$ by a compactly support isotopy so that
$\hat{{\CR}}^*(D^3_+ \setminus K)^{[i]}$ intersects
${\CR}^*(D^3_-\setminus K)^{[i]}$ transversally at a finite number of
intersection points. Denote the perturbed  diffeomorphism
by $\hat{f }_{\b }$. So its fixed points are all nondegenerated.

The Casson-Lin invariant of a knot $K = \overline{\b}$ is given by
counting the algebraic intersection number of
$\hat{{\CR}}^*(D^3_+ \setminus K)^{[i]}$ and
${\CR}^*(D^3_-\setminus K)^{[i]}$, or the algebraic number of
$\Fix(\hat{f }_{\b })$,
\[\lambda_{CL} (K) =  \lambda_{CL}(\b ) = Algebraic( \# \Fix(\hat{f }_{\b })) =
Algebraic( \# (
\hat{{\CR}}^*(D^3_+ \setminus K)^{[i]} \cap {\CR}^*(D^3_-\setminus K)^{[i]})) .\]

The results proved by Lin in \cite{lin} show that the Casson-Lin invariant 
$\lambda_{CL} (K) = \lambda_{CL}(\b )$
is independent of its braid representatives, i.e. $\lambda_{CL}(\b )$ is invariant
under the Markov moves of type I and type II on $\b$ and is one half of
the classical signature of the knot $K$.

\section{Nielsen number is a knot  invariant}

In this article we propose to count fixed points of  ${f }_{\b }$ in a Nielsen way  - using the classical  Nielsen numbers of ${f }_{\b }$.
 Nielsen counting of fixed points is a counting in the presence of the fundamental group.
In order to get an invariant of knots from braids, we have to verify that
Nielsen number
$ N(f_{\b})$  is invariant
under Markov moves. A Markov move of type I changes $\s \in B_n$ to
$\xi^{-1} \s \xi \in B_n$ for any $\xi \in B_n$, and the Markov move of
type II changes $\s \in B_n$ to $\s_n^{\pm } \s \in B_{n+1}$, or the inverses
of these operations. It is well-known that two braids $\b_1$ and $\b_2$ has
isotopic closure if and only if $\b_1$ can be changed to $\b_2$ by a sequence
of finitely many Markov moves \cite{bi}.
\begin{theorem} \label{invar}
If  $\overline{\b_1} = \overline{\b_2} = K$ as a knot, $\b_1 \in B_n, \b_2 \in
B_m$, then
\[ N(f_{\b_1})=N(f_{\b_2}).\]
So the Nielsen number 
$ N(f_{\b})$ is a knot invariant.
\end{theorem}
Proof: We only need to show that for $\b \in B_n$ with $\overline{\b}$ being
a knot $K$, the Markov moves of type I and type II on $\b$ provide either a cojugacy or a
isotopy of $f_{\b}$. Hence from the commutativity and  the invariance property under isotopy of the Nielsen numbers, we get that $ N(f_{\b})$
is an invariant of knot $K = \overline{\b}$.

Suppose we have the Markov move of type I: change $\b$ to $\xi^{-1} \b \xi$
for some $\xi \in B_n$. The element $\xi$ in $B_n$ induces a diffeomorphism
$\xi: Q_n \to Q_n$ is orientation preserving as observed by Lin in \cite{lin}.
Note that $B_n$ is generated by $\s_1, \cdots, \s_{n-1}$. For any $\s_i^{\pm }$,
the induced diffeomorphism $\s_i^{\pm } \times \s_i^{\pm }:
Q_n \times Q_n \to Q_n \times  Q_n$ is an orientation preserving  diffeomorphism. So $\xi $ is also a  orientation preserving
diffeomorphism since orientation  preserving
properties are invariant under the composition operation. Hence there is a
homeomorphism
\[\xi \times \xi : Q_n \times Q_n \to Q_n \times Q_n , \]
which commutes with the $SU(2)$-action and
\[ \xi \times \xi ({\CR}^*(S^2 \setminus K)^{[i]}) = {\CR}^*(S^2 \setminus K)^{[i]}
 \ \ \ (\mbox{changing variables by $\xi \times \xi$}),\]
\[\xi \times \xi ({\CR}^*(D_-^3 \setminus K)^{[i]}) =
{\CR}^*(D_-^3 \setminus K)^{[i]} \ \ \
(\mbox{in new coordinate $\xi(X_1), \cdots, \xi (X_n)$}),\]
\[ \xi \times \xi ({\CR}^*(D_+^3 \setminus K)^{[i]}) = 
{\CR}^*(D_+^3 \setminus K)^{[i]}\ \ \
(\mbox{in new coordinate $\xi(X_1), \cdots, \xi (X_n)$}),\]
as oriented manifolds. Let $g_{\xi} : {\CR}^*(S^2 \setminus K)^{[i]} \to
{\CR}^*(S^2 \setminus K)^{[i]}$ be the induced  homeomorphism, induced   from
$\xi \times \xi$ as coordinate changes . Hence we get a conjugacy relation
\[g_{\xi}^{-1}\circ f _{\b}\circ g_{\xi}  =  f_{\xi^{-1} \b \xi} , \]
from changing variables via $g_{\xi}$. Note that $Fix(f_{\xi^{-1} \b \xi})$
is identified with $Fix(f_{\b})$ under $g_{\xi}$. Thus the Markov move of type I preserves the conjugacy class of $ f_{\b}$
 Therefore by commutativity of the Nielsen number (see \cite{jb}) we have,
\begin{equation} \label{type1}
N(f_{\xi^{-1} \b \xi}) =
N(g_{\xi}^{-1}\circ f _{\b}\circ g_{\xi})
= N(f_{\b}) .
\end{equation}
It is clear that the argument goes through for the inverse operation of Markov
move of type I.

Suppose we have the Markov move of type II: change $\b$ to 
$\s_n \b \in B_{n+1}$. Recall that $\s_n (x_i) = x_i, 1 \leq i \leq n-1,
\s_n (x_n) = x_n x_{n+1} x_n^{-1}$ and $\s_n (x_{n+1}) = x_n$. We need to 
identify the Nielsen number  from the construction in
$\hat{H}_n$ into the one from $\hat{H}_{n+1}$. Following Lin \cite{lin}, there
is an imbedding $q: Q_n \times Q_n \to Q_{n+1} \times Q_{n+1}$ given by
\[q(X_1, \cdots, X_n, Y_1, \cdots, Y_n) = (X_1, \cdots, X_n,Y_n,
Y_1, \cdots, Y_n, Y_n) .\]
Such an imbedding commutes with the $SU(2)$-action and $q(H_n) \subset H_{n+1}$,
and induces an imbedding
\[ \hat{q}: \hat{H}_n (={\CR}^*(S^2 \setminus \overline{\b})^{[i]}) \to
\hat{H}_{n+1} (= {\CR}^*(S^2 \setminus \overline{\s_n \b})^{[i]}) .\]
Note that the symplectic structure of $\hat{H}_{n+1}$ restricted on
$\hat{q}(\hat{H}_n)$ is the symplectic structure on $\hat{H}_n$.
Under this imbedding, we have
$\hat{q}(f_{\b}): \hat{H}_{n+1} \to \hat{H}_{n+1}$ is given by
\begin{equation} \label{pb}
(X_1, \cdots, X_n, X_1, \cdots, X_n) \mapsto (X_1, \cdots, X_n, \b(X_n),
\b(X_1), \cdots, \b(X_n), \b(X_n)).
\end{equation}
The image of $\hat{q}(f_{\b})$ is invariant under the operation of $\s_n$.
Also the corresponding  diffeomorphism $f_{\s_n \b}$ is given by
\[f_{\s_n \b}(X_1, \cdots, X_n, X_{n+1}, X_1, \cdots, X_n, X_{n+1}) \]
\begin{equation} \label{pnb}
 = (X_1, \cdots, X_{n+1}, \b(X_1), \cdots, \b(X_{n-1}), \b(X_n) X_{n+1}
\b(X_n)^{-1}, \b(X_n) ) .
\end{equation}
Thus we have
\[ \hat{q}({\CR}^*(D_-^3  \setminus \overline{\b})^{[i]}) \subset
{\CR}^*(D_-^3  \setminus \overline{\s_n \b})^{[i]}, \ \
\hat{q}({\CR}^*(D_+^3  \setminus \overline{\b})^{[i]}) \subset
{\CR}^*(D_+^3  \setminus \overline{\s_n \b})^{[i]} . \]
The fixed points of $f_{\s_n \b}$ are elements
\[ \b(X_i) = X_i, 1 \leq i \leq n_1; \ \ \
 \b(X_n) X_{n+1} \b(X_n)^{-1} = X_n, \ \ \
\b(X_n) = X_{n+1} , \]
which is equivalent to $\b(X_i) = X_i, 1\leq i \leq n$, i.e.
\[ \Fix(f_{\s_n \b}) =\ Fix (\hat{q}(f_{\b})) = \Fix (f_{\b}) . \]

Then there is  a (Hamiltonian) isotopy
$\psi_t : \hat{H}_{n+1} (= (H_{n+1} \setminus S_{n+1})/SU(2)) \to \hat{H}_{n+1}$
between $\psi_{t_0} = \hat{q}(f_{\b})$ by (\ref{pb}) and
$\psi_{t_1} = f_{\s_n \b}$ (\ref{pnb})(see \cite{li,lin} for the explicit
constructions).  So we have
\begin{equation} \label{type2}
N(f_{\s_n \b}) = N(\hat{g}(f_{\b})) =
N(f_{\b}) .
\end{equation}
The first equality  is from the invariance property of Nielsen numbers under the  isotopy $\psi_t$ and the second
from the natural identification. We can similarly prove that
\[N(f_{\s_n^{-1}\b}) =
N(f_{\b}) .\]

\begin{rk}
It is known for a long time, that the problem of  computation of  Nielsen numbers is a  very difficult problem. By this reason, we strongly believe that the  Nielsen number $N(f_{\b})$ is a new knots invariant, which cannot be reduced to the  known knots invariants,  as it  happened   in the  case of the Casson-Lin invariant of knots.

\end{rk}

\begin{ex}
``Pillowcase''. For $n=2$ the irreducible variety  ${\CR}^*(S^2 \setminus K)^{[i]}$ is a 2-sphere with four cone points deleted(see \cite{lin}). So, in this case the space
${\CR}^*(S^2 \setminus K)^{[i]}$ is non-simply-connected  and the  Nielsen number $N(f_{\b})$ is   not trivial  for  general $\b$.

\end{ex}
\begin{que}
Does  the space ${\CR}^*(S^2 \setminus K)^{[i]}$ is simply-connected
 if  $n>2$?
\end{que}
The author was informed by Hans Boden  that the theorem 8.3 in \cite{furste} about  a moduli space of stable parabolic bundles over 2- sphere
 with marked points with given degree  and weights may be very usefull for the full   answer
on  this question.

\end{document}